\documentclass[11pt]{article}
\usepackage{latexsym}
\usepackage{amsfonts}
\usepackage[all]{xy}
\usepackage{amssymb, mathrsfs, amsfonts, amsmath}
\usepackage{amsbsy}
\newtheorem{thm}{Theorem}
\newtheorem{lem}[thm]{Lemma}

\begin{document}

\title{Essential spectrum of a class of Riemannian manifolds}

\author{Luiz Ant\^onio C. Monte \and Jos\'e Fabio B. Montenegro}

\maketitle

\noindent{\bf Abstract } In this paper we consider a family of Riemannian manifolds, not necessarily complete, with curvature conditions in a neighborhood of a ray. 
Under these conditions we obtain that the essential spectrum of the Laplace operator contains an interval. The results presented in this paper allow to determine 
the spectrum of the Laplace operator on unlimited regions of space forms, such as horoball in hyperbolic space and cones in Euclidean space.
\vspace{0.3cm}

\vspace{0.3cm}


\section{Introduction}
Let $M$ be a simply connected Riemannian manifold. The Laplace operator $\Delta :C_{0}^{\infty}(M)\to C_{0}^{\infty}(M)$, defined as
$\Delta u=\mathrm{div}(\mathrm{grad}\,u)$, is a second order elliptic operator and it has a unique extension $\Delta$ to an unbounded self-adjoint operator
on $L^2(M)$. Since $-\Delta$ is positive and symmetric, its spectrum is the set of $\lambda \geq0$ such that
$\Delta +\lambda I$ does not have bounded inverse. Sometimes we say spectrum of M rather than spectrum of $-\Delta$. 
One defines the essential spectrum $\sigma_{ess}(-\Delta )$ to be
those $\lambda$ in the spectrum which are either accumulation points of the spectrum or eigenvalues of infinite multiplicity.

It is well-known that if $M$ is an $n$-dimensional simply connected complete manifold with constant curvature $-c\leq 0$,
then  its essential spectrum coincides with the spectrum, being such the interval
$[(n-1)^2c/4, \infty )$. Moreover, the Decomposition Principle \cite{DoLi} says that the essential spectrum is invariant under compact pertubations of the metric 
on M and is thus a function of the geometry of the ends.
Therefore, it becomes natural the search of geometric conditions, of the ends of the surface, that will determine the essential spectrum of the Laplacian.
In 1981, Harold Donnelly in \cite{Donnelly}
studied the essential spectrum of manifolds which curvature approaches a constant $-c\leq 0$ at infinity. It was shown that the essential spectrum
is $[(n-1)^2c/4, \infty )$ if either (i) $M$ is simply connected and negatively curved or (ii) $M$ is a surface with finitely generated
fundamental group and an additional decay condition is satisfied for $K+c \to 0$, where $K$ is the Gaussian curvature. In 1992,
Escobar and Freire \cite{Escobar} proved that the spectrum of the Laplacian is $[0,\infty)$, using that the sectional curvature is non-negative
and the manifold satisfies some additional conditions. In \cite{Zhou}, Detang Zhou proved
that those additional conditions could be removed.
 In 1994, Li \cite{Li}  proved
$\sigma_{ess}(-\Delta) = [0, \infty)$ if $M$ has nonnegative Ricci curvatures and a pole.
Chen and Zhiqin Lu \cite{Chen} proved the same result when the radial sectional curvature is non-negative. Among other results, in \cite{Donn}, Donnelly proved that the
essential spectrum is $[0,\infty)$ for manifold with non-negative Ricci curvature and Euclidean volume growth. In 1997, Kumura \cite[Theorem 1.2]{Kumura} presented the 
following result: if $r$ is the distance function from a pole,
then  $\sigma_{ess}(-\Delta) = [c^{2}/4, \infty)$ provided
\begin{equation}
 \label{Deltar}
 \displaystyle\lim_{n \rightarrow \infty}\displaystyle\sup_{r\geq n}|\Delta r - c | = 0.
\end{equation}
Kumura also shows that this result recovers almost all the previous ones mentioned above.
In 1997, J. Wang in \cite{wang} proved that, if the Ricci curvature of a manifold $M$ satisfies $Ric(M)\geq -\delta / r^2$, where $r$ is the distance to a fixed point, 
and $\delta$ is a positive number depending only on the dimension, then the $L^p$ essential spectrum of $M$ is $[0,\infty )$ for any $p\in [1,+\infty ]$. 
In 2011, Zhiqin Lu and Detang Zhou in \cite{Zhiqin and Zhou} proved that the $L^{p}$ essential spectrum of the Laplacian is
$[0,+\infty)$ on a noncompact complete Riemannian manifold when
$$
 \displaystyle\liminf_{x \rightarrow \infty}Ric_{M}(x) = 0.
$$
 
The last result we want to recall is a theorem of Donnelly and Li. Fix a point $p\in M$ and write $\bar{K}(r)=\sup \{K(x,\pi )|d(p,x)\geq r\}$ where $K(x,\pi )$ is
the sectional curvature of a two plane $\pi$ in $T_xM$. Then the following theorem was proved in \cite{DoLi}.

\begin{thm}[H. Donnelly and P. Li] 
\label{thmDon}
Let $M$ be a complete Riemannian manifold and suppose $\bar{K}(r)\to -\infty$ as $r\to \infty$. 
Then the essential spectrum of $\Delta$ is empty provided one of the following 
two side conditions is satisfied: (i) $M$ is simply connected and negatively curved. (ii) $M$ is two dimensional and the fundamental group of $M$ is finitely generated.
\end{thm}

In this paper we consider a family of Riemannian manifolds, not necessarily complete, with curvature conditions like (\ref{Deltar}), 
but not uniform fashion over M. We require that it holds just in a neighborhood of a ray. 
With these conditions we obtain that the essential spectrum of the Laplacian contains an interval. This way, we are able to construct an
example of a two dimensional negatively curved Riemannian manifold satisfying
$$
\lim_{r\to \infty}K(r,\theta )(\frac{\partial}{\partial r}, \frac{\partial}{\partial \theta})=-\infty 
$$
for all $\theta \neq 0$ and such that the essential spectrum of this manifold contains the interval $[1/4, \infty)$, see the appendix. This example indicates that, 
in order for the essential spectrum to be empty, some type of geometric global conditions like Theorem \ref{thmDon}
should be necessary. The following Theorem provides some geometric conditions just on a neighborhood of a ray in order to guarantee that the 
essential spectrum is non-trivial.
\begin{thm}
\label{thm1}
Let $M$ be an n-dimensional Riemannian manifold. Suppose that, in geodesic spherical coordinates,
its metric can be written as $$g_M = dr^{2} + \psi^{2}(rw)g_{\mathbb{S}^{n-1}}$$ on $ C_a(N) = \left\{rw\; ;\,\, dist_{\mathbb{S}^{n-1}}(w,N) < c_2\,e^{-ar} \right\}$, 
where $g_{\mathbb{S}^{n-1}}=ds^2+\sin^2\!\!s \; g_{\mathbb{S}^{n-2}}$ and $s$
is the distance function in $\mathbb{S}^{n-1}$ to $N \in\mathbb{S}^{n-1}$. Furthermore, the function $\psi$ satisfies 
\begin{description}
\item[i)] $\displaystyle\lim_{{r \to \infty}\atop{w \rightarrow N}}\displaystyle\frac{\psi_{r}(rw)}{\psi(rw)}= c > 0$ uniform on the set $C_a(N)$, and $c>a$
\item[ii)] $\left|\displaystyle\frac{\psi_{s}(rw)}{\psi(rw)}\right| \leqslant c_1$ for all $rw \in C_a(N)$
\end{description}

\noindent Then $[(n-1)^{2}c^{2}/4, \infty) \subset \sigma_{ess}(-\Delta )$.
\end{thm}

We observe that changing the metric, outside the set $ C_a(N)$, or changing the topology of $M$,  the interval $[(n-1)^{2}c^{2}/4, \infty)$ 
remains contained in the essential spectrum of $M$. In fact, we prove that $[(n-1)^{2}c^{2}/4, \infty)$ is contained in the spectrum
of $C_a(N)$. As an application of this theorem we can show that horoball and hyperbolic space has the same spectrum.
For this purpose it is sufficient to show that horoball
$\left\{rw ;\,\, \sin^2\!\!s\, e^r \leq (1 + \cos s)^2\right\} $ 
 contains a set of the form  $\left\{rw ;\,\,\mathrm{ dist}_{\mathbb{S}^{n-1}}(w,N)= s \leq e^{-r/2}\right\}$.

In the case where $c=0$, the following theorem implies that the essential spectrum of a cone 
$\left\{rw\, ;\,\, dist_{\mathbb{S}^{n-1}}(w,N) < c_2 \right\}$ in $\mathbb{R}^n$  is the interval $[0,\infty )$.

\begin{thm}
\label{thm3}
Let $M$ be a n-dimensional Riemannian manifold. Suppose that, in geodesic spherical coordinates, 
its metric can be written as $$g_M = dr^{2} + \psi^{2}(rw)g_{\mathbb{S}^{n-1}}$$ on the set $ C_0(N) = \left\{\,\,rw ;\,\, dist_{\mathbb{S}^{n-1}}(w,N) < c_2 \right\}$, 
where $g_{\mathbb{S}^{n-1}}=ds^2+\sin^2\!\!s \; g_{\mathbb{S}^{n-2}}$ and $s$
is the distance function in $\mathbb{S}^{n-1}$ to $N \in\mathbb{S}^{n-1}$. Furthermore, assume that function $\psi$ satisfies 
\begin{description}
\item[i)] $\displaystyle\lim_{r \to \infty}\displaystyle\frac{\psi_{r}(rw)}{\psi(rw)}= 0$ to each $w$ such that $\mathrm{ dist}_{\mathbb{S}^{n-1}}(w,N) < c_2$
\item[ii)] $\displaystyle\lim_{r \to \infty}\psi(rw)= +\infty$ to each $w$ such that $\mathrm{ dist}_{\mathbb{S}^{n-1}}(w,N) < c_2$
\item[iii)] $\left|\displaystyle\frac{\psi_{s}(rw)}{\psi(rw)}\right| \leqslant \displaystyle\frac{c_1}{r^\gamma}$ for all $rw \in C_0(N)$ for some $\gamma > 1$.

\end{description}
\vspace{0.3cm}
\noindent Then $\sigma_{ess}(-\Delta )=[0,\infty )$.
\end{thm}

In both theorems, the function $\psi(rw)$ on the metric satisfies the following conditions:
$$\psi(0) = 0 \,\,\mbox{and}\,\, \psi(rw) > 0. $$
For more detail, see reference \cite{Chen}. The average curvature of the geodesic sphere of $M$ is given by
$$
\displaystyle\frac{\psi_{r}(rw)}{\psi(rw)}=\frac{1}{n-1}\;\Delta r.
$$

\section{Spectral Theory}

A linear operator on a Hilbert space $\mathcal{H}$ is a pair consisting of a dense linear subspace $\mathrm{Dom}(A)$ of $\mathcal{H}$ together with
a linear map $A:\mathrm{Dom}(A)\to \mathcal{H}$. The adjoint operator $A^*$ is determined by the condition that $\langle Au,v\rangle =\langle u,A^*v\rangle$
for all $u\in \mathrm{Dom}(A)$ and $v\in \mathrm{Dom}(A^*)$. The domain of $A^*$ is defined to be the set of all $v$ for which there exists $w\in \mathcal{H}$ such that
$\langle Au,v\rangle =\langle u,w\rangle$, for all $u\in \mathrm{Dom}(A)$. We say that $A$ is self-adjoint if $A=A^*$. The spectrum of a linear operator $A$,
$\sigma(A)$, is defined as follows. We say that a complex number $z$ does not lie in $\sigma(A)$ if the operator $(z-A)$ maps $\mathrm{Dom}(A)$ one-one onto $\mathcal{H}$,
and the inverse $(z-A)^{-1}$ is bounded. The spectrum of any self-adjoint operator is real and non-empty. A complex number is said to be an eigenvalue of such an
operator $A$ if there exists a non-zero $u\in \mathrm{Dom}(A)$ such that $Au=\lambda \,u$. It is entirely possible that no point of the spectrum of $A$ is an eigenvalue.
The discrete spectrum $\sigma_d(A)$ is defined as the set of all eigenvalues $\lambda$ of finite multiplicity which are isolated point of the spectrum. The essential spectrum is the set
$\sigma_{ess}(A)=\sigma(A)/\sigma_d(A)$. A characterization of the essential spectrum is given in following lemma which is a consequence of the spectral theorem 
\cite[Lemma 8.4.1, p.167]{Davies}.
\begin{lem}
\label{lemma}
 Let $A$ be a self-adjoint operator acting on the Hilbert space $\mathcal{H}$ and let $\lambda \in \mathbb{R}$. The following are equivalent:
\begin{description}
 \item[i)] $\lambda \in \sigma_{ess}(A)$
 \item[ii)] For all $\epsilon >0$ there exists a subspace $L_\epsilon \subset \mathrm{Dom}(A)$ with $dim (L_\epsilon )=\infty$ and such that 
$\| Au -\lambda u \|\leq \epsilon \, \|u\| $ for all $u\in L_\epsilon$.
 \end{description}
\end{lem}

\section{Proof of Theorem \ref{thm1}}
First we study the behavior of the function $\psi$ on the set $C_a(N)$. Let us prove that for any $\eta >0$ there are $r_{\eta} > 0$ such that
\begin{equation}
\label{est1psi}
C_{1}\,e^{(c - \eta)r}\leqslant \psi(rw)\leqslant C_{2}\,e^{ (c + \eta)r} 
\end{equation}
and
\begin{equation}
\label{est-quopsi}
\displaystyle\frac{1}{2} \leqslant \displaystyle\frac{\psi(rw)}{\psi(rN)} \leqslant \displaystyle\frac{3}{2} 
\end{equation}
for all $r \geq r_{\eta}$ and $rw \in C_a(N)$, where $C_{1}$ and $C_{2}$ are positive constants.
In fact, by the limit in the item i) of the Theorem \ref{thm1}, for any $ \eta > 0$, there is $r_{0}$ such that
$$
c - \eta \leqslant \displaystyle\frac{\psi_{r}(rw)}{\psi(rw)}\leqslant c + \eta
$$
for all $r \geq r_{0}$ and $rw \in C_a(N)$. Integrating the inequality above from $r_ {0}$ to $r$, we obtain
\begin{equation}
\label{estquopsipsi}
e^{(c - \eta)(r-r_{0})}\leqslant \displaystyle\frac{\psi(rw)}{ \psi(r_{0}w)}\leqslant e^{ (c + \eta)(r-r_{0})}
\end{equation}
for all $r \geq r_{0}$ and $rw \in C_a(N)$. By continuity and positivity of the function $w \mapsto \psi(r_{0}w)$,
$$\displaystyle\inf_{w \in\mathbb{S}^{n-1}} \psi(r_{0}w) > 0$$
and by (\ref{estquopsipsi})
$$
0 < C_{1}\,e^{(c - \eta)r}\leqslant \psi(rw)\leqslant C_{2}\,e^{ (c + \eta)r}
$$
for all $r\geqslant r_{0}$ and $rw\in C_a(N)$, where 
$$C_{1} = \displaystyle\inf_{w \in \mathbb{S}^{n-1}} \psi(r_{0}w)e^{ -(c - \eta)r_{0}}\;\;\; \mathrm{and}\;\;\;
C_{2} = \displaystyle\sup_{w \in \mathbb{S}^{n-1}} \psi(r_{0}w)e^{ -(c + \eta)r_{0}}.$$

To prove (\ref{est-quopsi}), consider $\alpha : [0,s] \rightarrow \mathbb{S}^{n-1}$ the geodesic such that $\alpha(0) = N$, $\alpha(s) = w$ and $\alpha'(t) = \partial/\partial s$. 
If $\lambda(t) = \psi(r\alpha(t))$, by Mean Value Theorem, there is $t_{0} \in (0,s)$ such that
$$ \lambda(s) - \lambda(0)= \lambda'(t_{0})\,s = s\, g_M(r\alpha(t_{0}))\left( \mathrm{grad }\psi,r\alpha'(t_{0}) \right)$$
$$\psi(rw) - \psi(rN) = s\,g_M (r\alpha(t_{0}))\left( \psi_{r}\displaystyle\frac{\partial}{\partial r}+\frac{\psi_{s}}{\psi^{2}}\,\displaystyle\frac{\partial}{\partial s} 
+\displaystyle\frac{1}{\psi^{2}}\, \mathrm{grad_{\mathbb{S}^{n-2}}}\, \psi,\; r\,\displaystyle\frac{\partial}{\partial s} \right)$$
$$\psi(rw) - \psi(rN) = r\,s\,\psi_{s}(r\alpha(t_{0})) $$
By the inequality in the item ii) of the Theorem (\ref{thm1}),
$$|\psi(rw) - \psi(rN)| \leqslant c_1\,r\,s\,\psi(r\alpha(t_{0}))$$
$$\left|\displaystyle\frac{\psi(rw)}{\psi(rN)} - 1\right|\leqslant\displaystyle\frac {c_1\,r\,s\,\psi(r\alpha(t_{0}))}{\psi(rN)}$$
Since $\mathrm{dist}_{\mathbb{S}^{n-1}}(\alpha(t_{0}),N)<s= \mathrm{dist}_{\mathbb{S}^{n-1}}(w,N) \leq c_{1}e^{-ar}$ and  using (\ref{est1psi}) 
$$\left|\displaystyle\frac{\psi(rw)}{\psi(rN)} - 1\right|\leqslant C \; \displaystyle\frac {r\,e^{(c+\eta)r}e^{-ar}}{e^{(c - \eta)r}} = 
C \;\displaystyle \frac{r}{e^{( a - 2\eta)r}} \longrightarrow 0 $$
when $r \rightarrow + \infty$, if  $0 < \eta < a/2 $.
Then there is $r_{\eta}\geqslant r_{0}$ which we obtain
$$
\displaystyle\frac{1}{2} \leqslant \displaystyle\frac{\psi(rw)}{\psi(rN)} \leqslant \displaystyle\frac{3}{2}
$$
for all $ r\geqslant r_{\eta}$ and $rw \in C_a(N)$.\\
In order to prove the Theorem \ref{thm1} we will construct, for any $\lambda > (n-1)^{2}c^{2}/4$ and $\epsilon > 0$, a sequence of functions
$(u_k)\subset C_0^{\infty}(M)$ with disjoint supports $\mathrm{supp} u_j\cap \mathrm{supp} u_k=\emptyset$, for all $j\neq k$, such that
\begin{equation}
\label{estLap-uk}
\|\Delta u_k +\lambda u_k\|_{L^2}\leq \epsilon \, \|u_k \|_{L^2}, \,\, k=1,2,\dots
\end{equation}
This implies that $[(n-1)^{2}c^{2}/4, \infty) \subset \sigma_{ess}(-\Delta )$ by the Lemma \ref{lemma} and the fact that $\sigma_{ess}(-\Delta )$ is  a
closed set. Indeed each function $u_k$ will have support on $C_a(N)$ and 
\begin{equation}
\label{def-uk}
u_k(rw)=f(r)g(s)
\end{equation}
where $s=\mathrm{dist}_{\mathbb{S}^{n-1}}(w,N)$. The function $f$ is defined by
\begin{equation}
\label{def-f}
f(r)=f(r,k,p)= F(r)h(r,k,p)
\end{equation}
where
\begin{equation}
\label{def-F}
F(r)=v(r)^{-1/2}\cos(\beta r),
\end{equation}
$\beta = \sqrt{\lambda - (n-1)^2c^2/4}$,
\begin{equation}
\label{def-v}
v(r)=\displaystyle\int_0^r \psi^{n-1}(\tau N)d\tau ,
\end{equation}
\begin{equation}
\label{def-h}
h(r)=h(r,k,p)=H\left(2(r-r_{k+2p})/(r_{k+4p}-r_{k})\right)
\end{equation}
the is a scaled cut-off function centered at  $r_{k+2p}$,
where $r_{k}=(2k+1)\pi/(2\beta )$ is the zero of the function $\cos(\beta r)$, 
and
$H \in C_{0}^{\infty}(\mathbb{R})$ is a cut-off function satisfying the conditions
\begin{center}
$\left\{
\begin{array}{l}
H\equiv 1 \;\mathrm{on}\; [-1/2,1/2]\\
H\equiv 0 \; \mathrm{on}\; \mathbb{R}\backslash [-1,1]\\
0\leqslant H\leqslant 1 \; \mathrm{on}\; \mathbb{R}.
\end{array}
\right.$
\end{center}
The function $g$ is defined by
\begin{equation}
\label{def-g}
g(s)=g(s,k,p)=H(s/\delta_{k,p})\cos (\pi s/\delta_{k,p})
\end{equation}
where $\delta_{k,p}=c_2\, e^{-ar_{k+4p}}.$
Now we will prove that there are $k_0>0$ and $p_0>0$ such that the functions $u_{k}=f(r,k,p)g(s,k,p)$ defined in (\ref{def-uk}) satisfies the inequality (\ref{estLap-uk}) for
all $k\geq k_0$ and $p=p(k)\geq p_0$.

The function $v(r)$ defined in (\ref{def-v}) satisfies
$$v'(r)= \psi^{n-1}(rN) \;\;\;\mathrm{and}\;\;\; v''(r)=(n-1)\psi^{n-2}(rN)\psi_{r}(rN).$$
 By 
(\ref{est1psi}), $\psi(rN)\geqslant M_{1}e^{(c-\eta)r} \,\,$ for all $\,\, r \geqslant r_{\eta}$ and $c> \eta >0$. Then
$$\lim_{r \to \infty}v'(r)= \lim_{r \to \infty}v(r) = +\infty $$
and
\begin{equation}
\label{limv'v}
\lim_{r \to \infty}\displaystyle\frac{v'(r)}{v(r)}= \lim_{r \to \infty}\displaystyle\frac{v''(r)}{v'(r)} =(n-1) \lim_{r \to \infty}\displaystyle\frac{\psi_r(rN)}{\psi(rN)} =(n-1)c.
\end{equation}
So there is a $r_{v}\geqslant r_{\eta}$ such that
\begin{equation}
\label{Est-v'/v}
 \displaystyle\frac{(n-1)c}{2} \leqslant \displaystyle\frac{v'(r)}{v(r)} \leqslant \displaystyle\frac{3(n-1)c}{2}\,\,\, ;\,\,\,\,\,\,\, \forall \;r\geqslant r_{v}.
 \end{equation}
\noindent The function $F(r)$ defined in (\ref{def-F}) satisfies
\begin{equation}
\label{LaplaceF}
\Delta  F + \lambda  F = A(r)F + B(rw)F'     
\end{equation}
where 
$$
A(r) = -\displaystyle\frac{1}{2} \displaystyle\frac{v''}{v'}\cdot \displaystyle\frac{v'}{v} + \displaystyle\frac{1}{4} \left(\displaystyle\frac{v'}{v}\right)^{2}+
\displaystyle\frac{(n-1)^2c^{2}}{4}
$$
and 
$$
B(rw)= (n-1)\displaystyle\frac{\psi_{r}}{\psi}(rw) - \displaystyle\frac{v'}{v}\;.
$$
Note that, by (\ref{limv'v}) and i) of the Theorem \ref{thm1},
\begin{equation}
\label{limitAandB}
\lim_{r \to \infty}A(r)=0 \;\;\;\mathrm{and} \;\;\; \lim_{{r \to \infty}\atop{w \rightarrow N}}B(rw) = 0
\end{equation}
uniform on $C_a(N)$.
\noindent By (\ref{LaplaceF}), the function $f$ defined in (\ref{def-f}) satisfies
\begin{equation}
\label{Laplace f}
\Delta f+\lambda  f= A(r)F\,h+ B(r,w)F'\,h+2 F'\,h'+F\,\Delta h  \; .    
\end{equation}
\noindent We have the following estimates for the function $h$ defined in (\ref{def-h}) 
\begin{equation}
\label{esth}
|h'|\leqslant \displaystyle\frac{\beta}{4 \pi p}sup|H'|\;\chi_h\;\;\;\;\;\;\;\;  \mbox{and}  \;\;\;\;\;\;\;\;  |h''|\leqslant \displaystyle\frac{\beta^2}{16 \pi^2 p^{2}}sup|H''|\;\chi_h.
\end{equation}
The Laplacian of the function $g=g(s)$, defined in (\ref{def-g}), is
\begin{equation}
\label{laplace g}
\Delta g = \displaystyle\frac{(n-3)\psi_{s}}{\psi^{3}}g'+ \displaystyle\frac{(n-2)\cot(s)}{\psi^{2}}g'+ \displaystyle\frac{1}{\psi^{2}}g''.
\end{equation}
We observed that
\begin{equation}
\label{est-g'}
|g'|\leqslant \displaystyle\frac{C}{\delta_{k,p}}\;\chi_{B(\delta_{k,p})}
\end{equation}
and
\begin{equation}
\label{est-g''}
|g''|\leqslant \displaystyle\frac{C}{\delta^{2}_{k,p}}\;\chi_{B(\delta_{k,p})}
\end{equation}
where $C$ is independent of $k$ and $p$, and $\chi_{B(\delta_{k,p})}:\mathbb{S}^{n-1}\to \mathbb{R}$ is the characteristic function of the set 
$B(\delta_{k,p})=\{w\in \mathbb{S}^{n-1};\;\mathrm{dist}_{\mathbb{S}^{n-1}}(w,N)\leq \delta_{k,p}\}$.

Finally, the function $u=u_k(rw) = f(r,k,p)\,g(s,k,p)$ satisfies
$$
\Delta u = (\Delta f )g+ f (\Delta g).
$$
Hence by (\ref{Laplace f}) and (\ref{laplace g}) we conclude
\begin{equation}
\label{Laplace u}
\Delta u + \lambda u =    
\end{equation}
$$
AFgh+ BF'\!gh+2 F'\!gh'+Fg\Delta h +\displaystyle\frac{(n-3)\psi_{s}}{\psi^{3}}fg'+ \displaystyle\frac{(n-2)\cot(s)}{\psi^{2}}fg'+ \displaystyle\frac{1}{\psi^{2}}fg'' 
$$

By (\ref{limitAandB}), given $\delta > 0$, there is $r_0>r_v$ such that $$|A(r)|\leqslant\delta \;\;\mathrm{and}\;\; |B(rw)|\leqslant\delta$$ 
for all $r\geqslant r_{0} $ and $rw \in C_a(N)$. 

By (\ref{esth}) there exists $r_{h}\geq r_{0}$ such that $$\|F'\,g\,h'\|_2\leqslant\delta \,\| \chi_h\,F'\,g\|_2\;\; \mathrm{and}\;\; \|F\,g\,\Delta h\|_2 \leqslant\delta \,\|\chi_h\,F\,g\|_2$$  
for all $r\geqslant r_{h}$, followed by (\ref{Laplace u}) 
\begin{equation}
\label{Delta u}
\|\Delta u + \lambda  u\|_2 \leqslant
\end{equation}
$$
 \delta \,\left(\|\chi_h\, F\,g\|_2+ \|\chi_h\,F'\,g\|_2 \right)+C\left\|\displaystyle\frac{\psi_{s}}{\psi^{3}}f\,g'\right\|_2 \!+ 
C\left \|\displaystyle\frac{\cot(s)}{\psi^{2}}f\,g'\right\|_2 \!+ \left\|\displaystyle\frac{1}{\psi^{2}}f\,g''\right\|_2
$$
for all  $r\geqslant r_{h}$ and $rw \in C_a(N)$.
We will need to use the technical lemma:
\begin{lem}
\label{lem5}
For the functions $F$, $f$, $g$ and $u$ defined previously, we have the following inequalities
\begin{description}
 \item[(a)] $\|\chi_h\,F\,g\|_2 \leqslant C\,\|u\|_2$
 \item[(b)] $\|\chi_h\,F'\,g\|_2 \leqslant C\,\|u\|_2$
 \item[(c)] $\left\|\displaystyle\frac{\psi_{s}}{\psi^{3}}f\,g'\right\|_2 \leqslant \displaystyle\frac{C}{\delta_{k,p}}\left[\displaystyle\inf_{C_{k,p}}|\psi|\right]^{-2}\|u\|_2$
 \item[(d)] $\left\|\displaystyle\frac{\cot(s)}{\psi^2}f\,g'\right\|_2 \leqslant \displaystyle\frac{C}{\delta_{k,p}^2}\left[\displaystyle\inf_{C_{k,p}}|\psi|\right]^{-2}\|u\|_2$
 \item[(e)] $\left\|\displaystyle\frac{1}{\psi^2}fg''\right\|_2 \leqslant \displaystyle\frac{C}{\delta_{k,p}^2}\left[\displaystyle\inf_{C_{k,p}}|\psi|\right]^{-2}\|u\|_2$
 \end{description}
where $C_{k,p}=\{rw\,;\;r_k\leq r \leq r_{k+4p},\; w\in B(\delta_{k,p})\}$ and $C$ is a positive constant independent of $k$ and $p$.
\end{lem}
\textbf{Proof of Lemma}: Observe that
$$
\|u\|_2^2=\int_{B(\delta_{k,p})}g^{2}(s)\int_{r_{k}}^{ r_{k+4p}}\!\! \cos^2(\beta r)\, h^2(r)\displaystyle\frac{\psi^{n-1}(rw)}{v(r)}\,dr\,dw \geqslant
$$
$$
\geqslant  \int_{B(\delta_{k,p})}g^{2}(s)\int_{r_{k+p}}^{ r_{k+3p}}\!\! \cos^2(\beta r)\, \displaystyle\frac{\psi^{n-1}(rw)}{\psi^{n-1}(rN)}\cdot\displaystyle\frac{v'(r)}{v(r)}\,dr\,dw
$$
and
$$
\|\chi_hF\,g\|_2^2= \int_{B(\delta_{k,p})}g^{2}(s)\int_{r_{k}}^{ r_{k+4p}} \cos^2(\beta r)\,\displaystyle\frac{\psi^{n-1}(rw)}{v(r)}\,dr\,dw.
$$
We observe of the estimates (\ref{Est-v'/v}) and (\ref{est-quopsi}) that
\begin{equation}
\label{estugcos}
\|u\|_2^2 \geqslant \displaystyle\frac{(n-1)c}{2} \left(\displaystyle\frac{1}{2}\right)^{n-1}\!\! \int_{B(\delta_{k,p})}g^{2}(s)\int_{r_{k+p}}^{ r_{k+3p}}\!\! \cos^2(\beta r)\,dr\,dw
\end{equation}
and
$$
\|\chi_hF\,g\|_2^2 \leqslant \displaystyle\frac{3(n-1)c}{2} \left(\displaystyle\frac{3}{2}\right)^{n-1}\!\! \int_{B(\delta_{k,p})}g^{2}(s)\int_{r_{k}}^{ r_{k+4p}}\!\! \cos^2(\beta r)\,dr\,dw.
$$
Moreover,
$$
\int_{r_{k}}^{ r_{k+4p}}\!\cos^2(\beta r)\, dr =2\int_{r_{k+p}}^{ r_{k+3p}}\!\cos^2(\beta r)\, dr.
$$
The last two inequalities imply
\begin{equation}
\label{des1lem}
\|\chi_hF\,g\|_2\leqslant \sqrt{2}\;3^{n/2}\,\,\|u\|_2.
\end{equation}
For the second inequality, using integration by parts
$$
\int_{r_{k}}^{ r_{k+4p}}\!\!F'(r)^{2}\, \psi^{n-1}(rw)\,dr =-\int_{r_{k}}^{ r_{k+4p}}\!\!F(r)\,\Delta F(r)\,\psi^{n-1}(rw)\,dr
$$
through of the equality (\ref{LaplaceF}) 
$$
\|\chi_h\,F'\,g\|_2^2 =\!\int_{B(\delta_{k,p})}\!\!\!g^{2}(s)\int_{r_{k}}^{ r_{k+4p}}\!\!\!\!F\,[\lambda F\!-\!A(r)F\!-\!B(r,w)F'\,]\;\psi^{n-1}(rw)dr\,dw
$$
of the limits in (\ref{limitAandB}) it follows
$$
\|\chi_h\,F'\,g\|_2^2\leqslant (\lambda +1) \int_{B(\delta_{k,p})}\!\!\!g^{2}(s)\int_{r_{k}}^{ r_{k+4p}}\!\!\! F^{2}\;\psi^{n-1}(rw)dr\,dw
$$
$$
+ \int_{B(\delta_{k,p})}\!\!\!g^{2}(s)\int_{r_{k}}^{ r_{k+4p}}\!\!\!|F|\, |F'|\psi^{n-1}(rw)dr\,dw
\leqslant (\lambda +\displaystyle\frac{3}{2}) \|\chi_h F\,g\|_2^2 +\displaystyle\frac{1}{2} \|\chi_h\,F'\,g\|_2^2
$$
$$
\|\chi_h\,F'\,g\|_2^2\leqslant
(2\lambda +3)\,\,\|\chi_h\,F\,g \|_2^2
$$
and of the inequality (\ref{des1lem}) we obtain
$$ \|\chi_h\,F'\,g\|_2 \leqslant (2\lambda +3)\sqrt{2}\,3^{n/2}\,\,\|u\|_2.$$
Now the third inequality
$$
\left\|\displaystyle\frac{\psi_{s}}{\psi^{3}}f\,g'\right\|^{2}_2 = 
\int_{B(\delta_{k,p})}\!\!\!|g'|^{2}\!\!\int_{r_{k}}^{ r_{k+4p}}\displaystyle\frac{\psi_{s}^{2}}{\psi^{6}}\, f^{2}(r)\,\psi^{n-1}\,dr\,dw
$$
By virtue of the estimate (\ref{est-g'}) and the hyphotesis ii) of the Theorem (\ref{thm1})
\begin{equation}
\label{estpsifg'}
\left\|\displaystyle\frac{\psi_{s}}{\psi^{3}}f\,g'\right\|^{2}_2 \leq \frac{C}{\delta^{2}_{k,p}}\left[\displaystyle\inf_{C_{k,p}}|\psi|\right]^{-4}
\int_{B(\delta_{k,p})}\!\int_{r_{k}}^{ r_{k+4p}}\!\!  f^2(r)\psi^{n-1}\,dr\,dw
\end{equation}

\noindent By definition, $f(r)=v^{-1/2}\cos(\beta r)h(r)$ and $v'(r)=\psi^{n-1}(rN)$, then
$$
\int_{B(\delta_{k,p})}\!\int_{r_{k}}^{ r_{k+4p}}\!\!  f^2(r)\psi^{n-1}\,dr\,dw =  
\int_{B(\delta_{k,p})}\!\int_{r_{k}}^{ r_{k+4p}}\!\! \cos^2(\beta r) h^2(r)\displaystyle\frac{\psi^{n-1}}{v}\,dr\,dw
$$
$$
\leq \int_{B(\delta_{k,p})}\!\int_{r_{k}}^{ r_{k+4p}}\!\! \cos^2(\beta r)\displaystyle\frac{\psi^{n-1}(rw)}{\psi^{n-1}(rN)}\frac{v'}{v}\,dr\,dw.
$$
By the estimates (\ref{Est-v'/v}) and (\ref{est-quopsi})
\begin{equation}
\label{estfcos}
\int_{B(\delta_{k,p})}\!\int_{r_{k}}^{ r_{k+4p}}\!\!  f^2(r)\psi^{n-1}\,dr\,dw \leq  C\int_{B(\delta_{k,p})}\!\int_{r_{k}}^{ r_{k+4p}}\!\! \cos^2(\beta r)\,dr\,dw
\end{equation}
where
$$
\int_{B(\delta_{k,p})}\!\,dw = \int_{\mathbb{S}^{n-2}}\int_{0}^{\delta_{k,p}}\!\sin^{n-2}s\,ds\,d\xi
$$
and $d\xi$ is the canonical measure of ${\mathbb{S}^{n-2}}$.

There exists $s_0>0$ such that 
$$
\displaystyle\frac{1}{2} \leqslant \displaystyle\frac{\sin s}{s}\leqslant \displaystyle\frac{3}{2}
$$
for all $0<s< s_{0}$. If $0< \delta_{k,p}< s_{0}$ we have
$$
\int_{B(\delta_{k,p})}\!\,dw \leq C \int_{\mathbb{S}^{n-2}}\int_{0}^{\delta_{k,p}}\!s^{n-2}ds\,d\xi =  
C\,\delta^{n-1}_{k,p}
$$
and also
$$
\int_{0}^{\delta_{k,p}/2}\cos^{2}\left(\displaystyle\frac{\pi s}{\delta_{k,p}}\right)\sin^{n-2}s \,ds 
\geqslant C\!\!\int_{0}^{\delta_{k,p}/2}\cos^{2}\left(\displaystyle\frac{\pi s}{\delta_{k,p}}\right)s^{n-2} \,ds $$
where
$$ \!\!\int_{0}^{\delta_{k,p}/2}\cos^{2}\left(\displaystyle\frac{\pi s}{\delta_{k,p}}\right)s^{n-2} \,ds =
\frac{\delta_{k,p}^ {n-1}}{\pi^{n-1}}\int_0^{\pi/2}\cos^2\!s\,s^{n-2}ds=C\delta^{n-1}_{k}.$$
Then we conclude
$$
\int_{B(\delta_{k,p})}\!\,dw 
\leq C \int_{\mathbb{S}^{n-2}}\int_{0}^{\delta_{k,p}/2}\cos^{2}\left(\displaystyle\frac{\pi s}{\delta_{k,p}}\right)\sin^{n-2}s \,ds\,d\xi
$$
by definition of the function $G_{\delta_{k,p}}$ follow that
$$
\int_{B(\delta_{k,p})}\!dw \leqslant C
\int_{\mathbb{S}^{n-2}}\int_{0}^{\delta_{k,p}/2}\!\!\!\!\!G^{2}_{\delta_{k,p}}(s)\cos^{2}\left(\displaystyle\frac{\pi s}{\delta_{k,p}}\right)\sin^{n-2}\!s\, ds\,d\xi
$$
$$
\leqslant C\!\!\int_{\mathbb{S}^{n-2}}\int_{0}^{\delta_{k,p}}\!\!G^{2}_{\delta_{k,p}}(s)\,\cos^{2}\left(
\displaystyle\frac{\pi s}{\delta_{k,p}}\right)\sin^{n-2}\!s \,ds\,d\xi = C\!\!\int_{B(\delta_{k,p})}\!\!g^{2}(s)\,dw.
$$
Using the last inequality in the the estimate (\ref{estfcos}) we have
\begin{equation}
\label{estfgcos}
\int_{B(\delta_{k,p})}\!\int_{r_{k}}^{ r_{k+4p}}\!\!  f^2(r)\psi^{n-1}\,dr\,dw \leq  C\int_{B(\delta_{k,p})}g^{2}(s)\!\int_{r_{k}}^{ r_{k+4p}}\!\! \cos^2(\beta r)\,dr\,dw
\end{equation}
using (\ref{estugcos}), we obten
\begin{equation}
\label{estfu}
\int_{B(\delta_{k,p})}\!\int_{r_{k}}^{ r_{k+4p}}\!\!  f^2(r)\psi^{n-1}\,dr\,dw \leq  C \|u\|^{2}_{2}
\end{equation}
by (\ref{estpsifg'}) verified
\begin{equation}
\label{estfg'u}
\left\|\displaystyle\frac{\psi_{s}}{\psi^{3}}f\,g'\right\|^{2}_2 \leq
\frac{C}{\delta^{2}_{k,p}}  \left[\displaystyle\inf_{C_{k,p}}|\psi|\right]^{-4} \|u \|^{2}_{2}
\end{equation}
Now we will show the fourth inequality, using a similar procedure like  the last inequality
$$
\left\| \frac{\cot s}{\psi^{2}}f\,g'\right\|^{2}_2 \leqslant 
\frac{C}{\delta^{2}_{k,p}} 
\left[\displaystyle\inf_{C_{k,p}}|\psi|\right]^{-4} \,
\int_{B(\delta_{k,p})}\!\int_{r_{k}}^{ r_{k+4p}}\!\!\cot^{2}\!s\cos^{2}(\beta r)\,dr\,dw.
$$
Since
$$
\int_{B(\delta_{k,p})}\cot^{2}\!s\;\,dw = \int_{\mathbb{S}^{n-2}}\!\!\sqrt{\mathrm{det}\xi}\int_{0}^{\delta_{k,p}}\!\!\sin^{n-4}\!s\; ds\,d\xi \leq 
$$
$$
 \leq C\int_{\mathbb{S}^{n-2}}\!\!\sqrt{\mathrm{det}\xi}\int_{0}^{\delta_{k,p}}\!\!s^{n-4} ds\,d\xi =
 C \delta^{n-3}_{k,p}\int_{\mathbb{S}^{n-2}}\!\!\!\!\sqrt{\mathrm{det}\xi}\,d\xi 
$$
and
$$
\int_{B(\delta_{k,p})}\!\!\!\!g^{2}(s)\,dw = 
\!\!\int_{\mathbb{S}^{n-2}}\sqrt{\mathrm{det}\xi}\int_{0}^{\delta_{k,p}}G^{2}_{\delta_{k,p}}(s)\cos^{2}\left(\displaystyle\frac{\pi s}{\delta_{k,p}}\right)\sin^{n-2}\!s\, ds\,d\xi 
$$
$$
\hspace{1,5cm}\geq\int_{\mathbb{S}^{n-2}}\sqrt{\mathrm{det}\xi}\int_{0}^{\delta_{k,p}/2}\cos^{2}\left(\displaystyle\frac{\pi s}{\delta_{k,p}}\right)\sin^{n-2}\!s\, ds\,d\xi 
$$
$$
\hspace{1,5cm}\geq C\int_{\mathbb{S}^{n-2}}\sqrt{\mathrm{det}\xi}\int_{0}^{\delta_{k,p}/2}\cos^{2}\left(\displaystyle\frac{\pi s}{\delta_{k,p}}\right)\!s^{n-2}ds\,d\xi
$$
notice the following
$$ \!\!\int_{0}^{\delta_{k,p}/2}\cos^{2}\left(\displaystyle\frac{\pi s}{\delta_{k,p}}\right)s^{n-2} \,ds =
\frac{\delta_{k,p}^ {n-1}}{\pi^{n-1}}\int_0^{\pi/2}\cos^2\!s\,s^{n-2}ds= C\delta^{n-1}_{k}.
$$
thus
$$
\int_{B(\delta_{k,p})}\!\!g^{2}(s)\,dw  \geq C \delta^{n-1}_{k}\!\int_{\mathbb{S}^{n-2}}\!\!\!\!\!\!\!\!\sqrt{\mathrm{det}\xi}\,d\xi.
$$
However
$$
\int_{B(\delta_{k,p})}\cot^{2}s\,ds\,dw \leq C \delta_{k,p}^{-2} \int_{B(\delta_{k,p})}\!\!g^{2}(s)\,dw 
$$
in this way we ensure
$$
\left\|\frac{\cot(s)}{\psi^{2}}f\,g'\right\|^{2}_2 \leqslant 
\frac{C}{\delta^{4}_{k,p}}\left[\displaystyle\inf_{C_{k,p}}|\psi|\right]^{-4} \,\int_{B(\delta_{k,p})}\!\int_{r_{k}}^{ r_{k+4p}}\!\!g^{2}(s)\cos^{2}(\beta r)\,dr\,dw 
$$
of inequality (\ref{estugcos}) we deduce
\begin{equation}
\label{estcotfg'u}
\left\|\frac{\cot(s)}{\psi^{2}}f\,g'\right\|^{2}_2 \leqslant \frac{C}{\delta^{4}_{k,p}}\left[\displaystyle\inf_{C_{k,p}}|\psi|\right]^{-4} \,\|u \|^{2}_{2} .
\end{equation}
To finish, using (\ref{est-g''}) and the same reasoning above
$$
\left\|\frac{1}{\psi^{2}}\,f\,g''\right\|_2^ 2  \leqslant 
\frac{C}{\delta^{4}_{k,p}}\left[\displaystyle\inf_{C_{k,p}}|\psi|\right]^{-4} \int_{B(\delta_{k,p})}\!\int_{r_{k}}^{ r_{k+4p}}\!\!f^{2}(r)\psi^{n-1}\,dr\,dw .
$$
By (\ref{estfu}), we conclude
\begin{equation}
\label{estfg"u}
\left\|\frac{1}{\psi^{2}}\,f\,g''\right\|_2^ 2  \leqslant \frac{C}{\delta^{4}_{k,p}}\left[\displaystyle\inf_{C_{k,p}}|\psi|\right]^{-4} \,\,\|u \|^{2}_{2} 
\end{equation}
Continuing the proof of the theorem, consider $p =\lfloor k/m \rfloor$( The party integer of $k/m$, where $m \in \mathbb{N}$), so by definition of $r_{k+4p}$ we obtain
$r_{k+4p} \leqslant (1 + 4/m)r_{k} + M$, jointly with the defintion $\delta_{k,p}$ and with the items c), d) and e) above the lemma \ref{lem5}
$$
\left\|\displaystyle\frac{(n-1)\psi_{s}}{\psi^{3}}f\,g'\right\|^{2}_2 \leqslant C\displaystyle\frac{e^{2a(1+4/m)r_{k} }}{e^{4(c-\eta)r_{k}}}\|u\|_{2}^ 2
$$
$$
\left\|\displaystyle\frac{(n-2)\cot(s)}{\psi^{2}}f\,g'\right\|^{2}_2 \leqslant C\displaystyle\frac{e^{4a(1+4/m)r_{k} }}{e^{4(c-\eta)r_{k}}}\|u\|_{2}^2
$$
$$
\left\|\displaystyle\frac{1}{\psi^{2}}f\,g''\right\|^{2}_2 \leqslant C\displaystyle\frac{e^{4a(1+4/m)r_{k} }}{e^{4(c-\eta)r_{k}}}\|u\|_{2}^2
$$
Given $\epsilon > 0$, there $\eta > 0$ and $m \in \mathbb{N}$ where $c-\eta > a(1+4/m)$ such that
$$
\left\|\displaystyle\frac{(n-1)\psi_{s}}{\psi^{3}}f\,g'\right\|_2,\left\|\displaystyle\frac{(n-2)\cot(s)}{\psi^{2}}f\,g'\right\|_2 ,\left\|\displaystyle\frac{1}{\psi^{2}}f\,g''\right\|_2  
\leq \epsilon\,\|u\|_{2}
$$
and jointly with the lemma \ref{lem5}, we deduce of (\ref{Delta u}) that
\begin{equation}
\label{estdeltau}
\|\Delta u+\lambda  u\|_2 \leq \epsilon\,\|u\|_{2}
\end{equation}

Consider the subspace spanned
\begin{center}
$G=[[u(k_0,p_0,.),u(k_0+4p_0,p_0,.),u(k_0+8p_0,p_0,.),\cdots]]$
\end{center}
where $supp( u(k_0+2^ip_0,p_0,.))\bigcap supp(u(k_0+2^jp_0,p_0,.))= \emptyset$, for $i \neq j $, so
$$
\|\Delta u +\lambda  u\|_2 < \epsilon \, \|u\|_2
$$
for all $u\in G$. By the lemma \ref{lemma}, $\lambda \in \sigma_{ess}(-\Delta )$, which concluded
\begin{equation}
\label{incspe}
[(n-1)^{2}c^{2}/4, \infty) \subseteq \sigma_{ess}(-\Delta )   
\end{equation}

\section{Proof of Theorem \ref{thm3}}
 Let us prove that there is $r_{0} > 0$ such that
\begin{equation}
\label{estquopsi}
\displaystyle\frac{1}{2} \leqslant \displaystyle\frac{\psi(rw)}{\psi(rN)} \leqslant \displaystyle\frac{3}{2} 
\end{equation}
for all $r \geq r_{0}$ and $rw \in C_0(N)$.
Consider $\mu : [0,s] \rightarrow \mathbb{S}^{n-1}$ the geodesic such that $\mu(0) = N$, $\mu(s) = w$ and $\mu'(t) = \partial/\partial s$.
If $\nu(t) = \psi(r\mu(t))$, by Mean Value Theorem, there is $t_{0} \in (0,s)$ such that
$$ \nu(s) - \nu(0)= \nu'(t_{0})\,s = s\, g_M(r\mu(t_{0}))\left( \mathrm{grad }\psi,r\mu'(t_{0}) \right)$$
$$\psi(rw) - \psi(rN) = s\,g_M (r\mu(t_{0}))\left( \psi_{r}\displaystyle\frac{\partial}{\partial r}+\frac{\psi_{s}}{\psi^{2}}\,\displaystyle\frac{\partial}{\partial s} 
+\displaystyle\frac{1}{\psi^{2}}\, \mathrm{grad_{\mathbb{S}^{n-2}}}\, \psi,\; r\,\displaystyle\frac{\partial}{\partial s} \right)$$
\begin{equation}
\label{psi(wN)}
\psi(rw) - \psi(rN) = r\,s\,\psi_{s}(r\mu(t_{0})) 
\end{equation}
Of the hypotese iii)
$$
\frac{\psi_{\tau}(r,\tau,\xi)}{\psi(r,\tau,\xi)} \leq \frac{c_1}{r^{\gamma}}
$$
in $C_0(N)$, integrating with respect to $\tau$ from 0 to s, we obten
$$
\psi(r,s,\xi)\leq \psi(r,0,\xi)e^{c_{1}s/r^\gamma} = \psi(rN)e^{c_{1}s/r^\gamma}
$$
thus of iii) and (\ref{psi(wN)})

$$|\psi(rw) - \psi(rN)| \leq \frac{c_{1}\,c_2\,e^{c_{1}c_2/r^\gamma}\psi(rN)}{r^{\gamma}}$$
since $\gamma > 1$, there $r_0$ such that
\begin{equation}
\label{estquopsi2}
\frac{1}{2}\leq \frac{\psi(rw)}{\psi(rN)} \leq \frac{3}{2}
\end{equation}
for all $r\geq r_0$ and $rw \in C_0(N)$.\\
For $\lambda > 0$, take $\alpha > 0$ such that $\lambda > \frac{(n-1)^2\alpha^2}{4}$. Consider the metric
$$
g_{M^{\alpha}} = dr^{2} + \psi^{\alpha}(rw)^2g_{\mathbb{S}^{n-1}}
$$
in $C_0(N)$, where $\psi^{\alpha}(rw) = e^{\alpha\,r }\psi(rw)$. By hypothesis i), the function $\psi^{\alpha}(rw)$ satisfies
$$\displaystyle\lim_{r \to \infty}\displaystyle\frac{\psi^{\alpha}_{r}(rw)}{\psi^{\alpha}(rw)}=\displaystyle\lim_{r \to \infty}\displaystyle\left(\frac{\psi_{r}(rw)}{\psi(rw)}+\alpha\right)=\alpha $$ to each $w$ such that $\mathrm{ dist}_{\mathbb{S}^{n-1}}(w,N) \leqslant c_2$. 

Given $\epsilon > 0$, construct a sequence $(u_k^{\alpha})$ in a manner analogous to $(\ref{def-uk})$ that satisfies (\ref{estLap-uk}), where
$f_k^{\alpha}$ is the same of $(\ref{def-f})$ with $v_{\alpha}(r)=\displaystyle\int_{0}^{r}\psi^{\alpha}(\tau N)^{n-1}d\tau$, $r_k = (2k+1)\pi/2\sqrt{\lambda}$ and $g(s)= H(s/c_2)\cos (\pi s/c_2)$, knowing that $g$ satisfies  (\ref{laplace g}) and $$supp\,g = B(c_2)=\left\{w\in \mathbb{S}^{n-1}\,\,;\,\,\mathrm{ dist}_{\mathbb{S}^{n-1}}(w,N) \leqslant c_2\right\} $$

Similarly to (\ref{Laplace u})
\begin{equation}
\label{iguallap2}
\Delta_{\alpha} u^{\alpha} + \lambda u^{\alpha} = A_{\alpha}(r)F^{\alpha}\,g\,h+ B_{\alpha}(rw)(F^{\alpha})'\,g\,h+2 (F^{\alpha})'\,g\,h'+F^{\alpha}\,g\,\Delta h +    
\end{equation}
$$
\hspace*{1.5cm} +\displaystyle\frac{(n-3)\psi^{\alpha}_{s}}{(\psi^{\alpha})^{3}}f^{\alpha}\,g'+ \displaystyle\frac{(n-2)\cot(s)}{(\psi^{\alpha})^{2}}f^{\alpha}\,g'+ \displaystyle\frac{1}{(\psi^{\alpha})^{2}}f^{\alpha}\,g''  \; .
$$
where
$$
A_{\alpha}(r) = -\displaystyle\frac{1}{2} \displaystyle\frac{v_{\alpha}''}{v_{\alpha}} + \displaystyle\frac{1}{4} \left(\displaystyle\frac{v_{\alpha}'}{v_{\alpha}}\right)^{2} + \frac{(n-1)^2\alpha^2}{4}
$$
and 
$$
B_{\alpha}(rw)= (n-1)\displaystyle\frac{\psi^{\alpha}_{r}}{\psi^{\alpha}}(rw) - \displaystyle\frac{v_{\alpha}'}{v_{\alpha}}\;.
$$
analogous to (\ref{limv'v})
\begin{equation}
\label{estquovalpha}
\lim_{r \to \infty}\displaystyle\frac{v_{\alpha}'(r)}{v_{\alpha}(r)}= \lim_{r \to \infty}\displaystyle\frac{v_{\alpha}''(r)}{v_{\alpha}'(r)} = (n-1)\alpha
\end{equation}
then
\begin{equation}
\label{limitAandB2}
\lim_{r \to \infty}A_{\alpha}(r)=0 \;\;\;\mathrm{and} \;\;\; \lim_{r \to \infty}B_{\alpha}(rw) = 0
\end{equation}
to each $w$ such that $\mathrm{ dist}_{\mathbb{S}^{n-1}}(w,N) \leqslant c_2$. 

Then, similarly the $(\ref{Delta u})$, given $\delta > 0$ there $r_0$ such that

\begin{equation}
\label{Deltaualpha2 }
\|\Delta_{\alpha} u^{\alpha} + \lambda  u^{\alpha}\|_{L^2(M^{\alpha})} \leqslant \delta \|\chi_h F^{\alpha}g\|_{L^2(M^{\alpha})}+ \delta \|\chi_h(F^{\alpha})'g\|_{L^2(M^{\alpha})} +
\end{equation}
$$
+\,C\left\|\displaystyle\frac{\psi^{\alpha}_{s}}{(\psi^{\alpha})^{3}}f^{\alpha}g'\right\|_{L^2(M^{\alpha})} + 
C\left \|\displaystyle\frac{\cot(s)}{(\psi^{\alpha})^{2}}f^{\alpha}g'\right\|_{L^2(M^{\alpha})} + \left\|\displaystyle\frac{1}{(\psi^{\alpha})^{2}}f^{\alpha}g''\right\|_{L^2(M^{\alpha})}
$$
for all $r\geqslant r_0$ and $rw \in C_0(N)$.
\begin{lem}
\label{lem6}
For the functions $F^{\alpha}$, $f^{\alpha}$, $g$ and $u^{\alpha}$ defined previously, we have the following inequalities
\begin{description}
 \item[(a)] $\|\chi_h\,F^{\alpha}\,g\|_2 \leqslant C\,\|u^{\alpha}\|_2$
 \item[(b)] $\|\chi_h\,(F')^{\alpha}\,g\|_2 \leqslant C\,\|u^{\alpha}\|_2$
 \item[(c)] $\left\|\displaystyle\frac{\psi^{\alpha}_{s}}{(\psi^{\alpha})^{3}}f^{\alpha}\,g'\right\|_2 \leqslant C\,r_k^{-\gamma}\left[\displaystyle\inf_{C_{k,p}}|\psi^{\alpha}|\right]^{-2}\|u^{\alpha}\|_2$
 \item[(d)] $\left\|\displaystyle\frac{\cot(s)}{(\psi^{\alpha})^2}f^{\alpha}\,g'\right\|_2 \leqslant C\left[\displaystyle\inf_{C_{k,p}}|\psi^{\alpha}|\right]^{-2}\|u^{\alpha}\|_2$
 \item[(e)] $\left\|\displaystyle\frac{1}{(\psi^{\alpha})^2}f^{\alpha}g''\right\|_2 \leqslant C\left[\displaystyle\inf_{C_{k,p}}|\psi^{\alpha}|\right]^{-2}\|u^{\alpha}\|_2$
 \end{description}
where $C_{k,p}=\{rw\,;\;r_k\leq r \leq r_{k+4p},\; w\in B(c_2)\}$ and $C$ is a positive constant independent of $k$ and $p$.
\end{lem}
\textbf{Proof of Lemma}: The items [(a)] and [(b)] is proved similarly to the items [(a)] and [(b)] of lemma \ref{lem5}. We prove item [(c)], the other follows similarly.

Now the third inequality
$$
\left\|\displaystyle\frac{\psi^{\alpha}_{s}}{(\psi^{\alpha})^{3}}f^{\alpha}\,g'\right\|^{2}_{L^2(M^{\alpha})} = 
\int_{B(c_2)}\!\!\!|g'|^{2}\!\!\int_{r_{k}}^{ r_{k+4p}}\displaystyle\frac{(\psi_{s}^{\alpha})^{2}}{(\psi^{\alpha})^{6}}\, f^{\alpha}(r)^2\,(\psi^{\alpha})^{n-1}\,dr\,dw
$$
By hyphotesis iii) and $\displaystyle\inf_{C_{k,p}}|\psi^{\alpha}| \leqslant |\psi^{\alpha}|$,  
$$
\left\|\displaystyle\frac{\psi^{\alpha}_{s}}{(\psi^{\alpha})^{3}}f^{\alpha}\,g'\right\|^{2}_{L^2(M^{\alpha})} \leq \frac{C}{r_k^{2\gamma}}\left[\displaystyle\inf_{C_{k,p}}|\psi^{\alpha}|\right]^{-4}
\int_{B(c_2)}\!\!\!|g'|^{2}\int_{r_{k}}^{ r_{k+4p}}\!\!  f^{\alpha}(r)^2(\psi^{\alpha})^{n-1}\,dr\,dw
$$
since that $f^{\alpha}(r)=v_{\alpha}^{-1/2}\cos(\beta r)h(r)$ and $v_{\alpha}'(r)=\psi^{\alpha}(rN)^{n-1}$, then
$$
\int_{B(c_2)}\!\!\!\!|g'|^{2}\int_{r_{k}}^{ r_{k+4p}}\!\!   f^{\alpha}(r)^2(\psi^{\alpha})^{n-1}\,dr\,dw =  
\int_{B(c_2)}\!\!\!\!|g'|^{2}\int_{r_{k}}^{ r_{k+4p}}\!\! \cos^2(\beta r) h^2(r)\displaystyle\frac{(\psi^{\alpha})^{n-1}}{v_{\alpha}}\,dr\,dw
$$
$$
\leq \int_{B(c_2)}\!\!\!\!|g'|^{2}\int_{r_{k}}^{ r_{k+4p}}\!\! \cos^2(\beta r)\displaystyle\frac{\psi^{\alpha}(rw)^{n-1}}{\psi^{\alpha}(rN)^{n-1}}\cdot\frac{v'_{\alpha}}{v_{\alpha}}\,dr\,dw.
$$
by the estimates (\ref{estquopsi2}), (\ref{estquovalpha}) and definition of $\psi^{\alpha}$ 
\begin{equation}
\label{estfg'cos2}
\left\|\displaystyle\frac{\psi^{\alpha}_{s}}{(\psi^{\alpha})^{3}}f^{\alpha}\,g'\right\|^{2}_{L^2(M^{\alpha})} \leq  \frac{C}{r_k^{2\gamma}}\left[\displaystyle\inf_{C_{k,p}}|\psi^{\alpha}|\right]^{-4}\!\!\!\int_{B(c_2)}\!\!\!|g'|^{2}\int_{r_{k}}^{ r_{k+4p}}\!\!\!\! \cos^2(\beta r)\,dr\,dw
\end{equation}
also have
\begin{equation}
\label{estucos2}
\|u^{\alpha}\|_{L^2(M^{\alpha})}^2 \geqslant C\,\alpha \!\! \int_{B(c_2)}g^{2}(s)\int_{r_{k+p}}^{ r_{k+3p}}\!\! \cos^2(\beta r)\,dr\,dw
\end{equation}
but
$$
\int_{B(c_2)}|g'(s)|^2\,dw \leqslant C\,\int_{B(c_2)}g(s)^2\,dw\,\,\mbox{and}\,\,\int_{r_{k}}^{ r_{k+4p}}\!\cos^2(\beta r)\, dr =2\int_{r_{k+p}}^{ r_{k+3p}}\!\cos^2(\beta r)\, dr.
$$
jontly with (\ref{estfg'cos2}) and (\ref{estucos2})
\begin{equation}
\label{estfg'u2}
\left\|\displaystyle\frac{\psi^{\alpha}_{s}}{(\psi^{\alpha})^{3}}f^{\alpha}\,g'\right\|^{2}_{L^2(M^{\alpha})} \leq  \frac{C}{r_k^{2\gamma}}\left[\displaystyle\inf_{C_{k,p}}|\psi^{\alpha}|\right]^{-4}\|u^{\alpha}\|_{L^2(M^{\alpha})}^2
\end{equation}

Continuing the proof of the theorem, in inequality (\ref{Deltaualpha2 }) we use the lemma \ref{lem6}, by  hypotheses ii), given $\epsilon > 0$ there $r_0$ such that
$$
\|\Delta_{\alpha} u^{\alpha}+\lambda  u^{\alpha}\|_{L^2(M^{\alpha})} \leq (\epsilon/M)\,\|u^{\alpha}\|_{L^2(M^{\alpha})}
$$
for all $r\geqslant r_0$ and $rw \in C_0(N)$.
Let us now define $u_k = \mu(r,\alpha)^{-1}u_k^{\alpha} \in C_{0}^{\infty}(M)$ where $\mu(r,\alpha)= e^{-(n-1)\alpha\,r/2}$, that notice $\|u_k \|_{L^2(M)} = \|u_k^{\alpha} \|_{L^2(M^{\alpha})}$. However
$$
\Delta_{\alpha} u_k^{\alpha}  =\left[ \Delta u_k +\left(\frac{1}{e^{2\alpha\,r}}-1\right)f_k\,\Delta g -(n-1)^2 \left(\frac{\alpha^2}{4} + \frac{\alpha}{2}\frac{\psi_r}{\psi}\right)u_k\right]\mu(r)
$$
where $f_k = \mu(r,\alpha)^{-1}f_k^{\alpha}$, and then
$$
\left\|\left[ \Delta u_k\!+\!\!\left(\frac{1}{e^{2\alpha\,r}}\!-\!1\right)f_k\,\Delta g -(n-1)^2 \!\!\left(\frac{\alpha^2}{4}\! + \!\frac{\alpha}{2}\frac{\psi_r}{\psi}\right)u_k\!+\!\lambda u_k\right]\mu(r)\right\|_{L^2(M^{\alpha})}\!\!\!\!\leq \!\!(\epsilon/M)\,\!\|u^{\alpha}\|_{L^2(M^{\alpha})}
$$
in other words
$$
\left\|\Delta u_k\!+\!\!\left(\frac{1}{e^{2\alpha\,r}}\!-\!1\right)f_k\,\Delta g -(n-1)^2 \!\!\left(\frac{\alpha^2}{4}\! + \!\frac{\alpha}{2}\frac{\psi_r}{\psi}\right)u_k+\!\lambda u_k\right\|_{L^2(M)} \!\!\leq (\epsilon/M)\,\|u\|_{L^2(M)}
$$
using triangular inequality
$$
\|\Delta u_k\!+\!\lambda\,u_k\|_{L^2(M)} \!\!\leq (\epsilon/M)\,\|u\|_{L^2(M)}+C\|f_k\,\Delta g \|_{L^2(M)}+C\left\|\left(\frac{\alpha^2}{4}\! + \!\frac{\alpha}{2}\frac{\psi_r}{\psi}\right)u_k\right\|_{L^2(M)}
$$
Note that
$$
\|f_k\,\Delta g \|_{L^2(M)} = \|f^{\alpha}_k\,\Delta g \|_{L^2(M^{\alpha})} =
$$
$$
\int_{B(c_2)}\!\int_{r_{k}}^{ r_{k+4p}}\!\!   f^{\alpha}(r)^2\Delta g^2(\psi^{\alpha})^{n-1}\,dr\,dw =  
\int_{B(c_2)}\!\int_{r_{k}}^{ r_{k+4p}}\!\! \cos^2(\beta r) h^2(r)\Delta g^2\displaystyle\frac{(\psi^{\alpha})^{n-1}}{v_{\alpha}}\,dr\,dw
$$
$$
\leq \int_{B(c_2)}\!\int_{r_{k}}^{ r_{k+4p}}\!\! \cos^2(\beta r)\Delta g^2\displaystyle\frac{\psi^{\alpha}(rw)^{n-1}}{\psi^{\alpha}(rN)^{n-1}}\cdot\frac{v'_{\alpha}}{v_{\alpha}}\,dr\,dw.
$$
by the estimates (\ref{estquopsi2}), (\ref{estquovalpha}) and definition of $\psi^{\alpha}$ 
$$
\|f_k\,\Delta g \|_{L^2(M)}\leq C\alpha\int_{B(c_2)}\!\int_{r_{k}}^{ r_{k+4p}}\!\! \cos^2(\beta r)\Delta g^2\,dr\,dw 
$$
for (\ref{laplace g}), hypothesis iii), $r\geqslant r_0$ and  $supp\,g = supp\,g' = supp\,g'' = B(c_2)$, it follows that
$$
\int_{B(c_2)}\!\Delta g^2\,dw \leqslant\int_{B(c_2)}\!\left[\left|\displaystyle\frac{(n-3)\psi_{s}}{\psi^{3}}\right||g'|+ \displaystyle\frac{(n-2)|cot(s)|}{\psi^{2}}|g'|+ \displaystyle\frac{1}{\psi^{2}}|g''|
\right]^2\,dw \leqslant 
$$
$$
\leqslant\!\! \frac{C}{(\displaystyle\inf_{C_{k,p}}|\psi|)^4}\!\!\int_{B(c_2)}\!\left[\left|\displaystyle\frac{(n-3)c_1}{r_0^{\gamma}}\right||g'|+ (n-2)|cot(s)||g'|+ |g''|
\right]^2\!\!dw \leqslant \!\frac{C}{(\displaystyle\inf_{C_{k,p}}|\psi|)^4}\int_{B(c_2)}\!\! g^2\!dw
$$
thus
$$
\|f_k\,\Delta g \|_{L^2(M)}\leqslant \!\frac{C\alpha}{(\displaystyle\inf_{C_{k,p}}|\psi|)^4}\int_{B(c_2)}g(s)^2\!\int_{r_{k}}^{ r_{k+4p}}\!\! \cos^2(\beta r)\,dr\,dw 
$$
through inequality (\ref{estucos2})
$$
\|f_k\,\Delta g \|_{L^2(M)}\leqslant \!\frac{C}{(\displaystyle\inf_{C_{k,p}}|\psi|)^4}\|u_k\|_{L^2(M)}
$$
of the hypothesis ii), we have
$$
\|f_k\,\Delta g \|_{L^2(M)}\leqslant C\epsilon\|u_k \|_{L^2(M)} 
$$
doing $\alpha ,|\psi_{r}/\psi| < \epsilon$, we conclude
\begin{equation}
\|\Delta u_k +\lambda u_k\|_{L^2(M)}\leq \epsilon \, \|u_k \|_{L^2(M)}, \,\, k=1,2,\dots
\end{equation}
by the lemma \ref{lemma}, that guarantee
\begin{equation}
\label{incspe2}
\sigma_{ess}(-\Delta ) = [0, \infty)
\end{equation}

\section{Appendix}
Consider $\mathbb{R}^{2} = \left\{(r\cos\theta, r\sin \theta)\,\,; \,\, r \geq 0 \,\,, \,\,\theta \in [0,2\pi]\right\}$ with metric given by 
$g = dr^{2}+\psi^{2}(r,\theta)g_{\mathbb{S}^{1}}$, since $\psi(r,\theta) =r e^{r^{2}g(\theta) + r}$ where $g(\theta) = \sin^{2}(\theta/2)$. Calculating the derivatives of $\psi$
$$
\psi_{r} = \frac{\psi}{r} + (2rg(\theta) + 1)\psi\,\,\,\,\mbox{and} \,\,\,\, \psi_{rr} = \frac{\psi_{r}}{r}-\frac{\psi}{r^{2}}+ 2g(\theta)\psi + (2rg(\theta) + 1)\psi_{r}
$$
then we deduce
$$
\left|\frac{\psi_{r}}{\psi} - 1\right| = \frac{1}{r} + 2r(\theta/2)^{2}\frac{\sin^{2}(\theta/2)}{(\theta/2)^{2}} 
$$
and
$$
K(r,\theta) = -\frac{\psi_{rr}}{\psi}(r,\theta) = -6g(\theta) -\frac{2}{r} - (2rg(\theta) + 1)^{2} 
$$
In the set $C(a)= \left\{ (r\cos\theta , r\sin\theta)\,\,;\,\, |\theta|\leq e^{-ar} \right\}$ we have that
$$
\left|\frac{\psi_{r}}{\psi} - 1\right| \leq Cre^{-2ar} \rightarrow 0
$$
when $r \rightarrow +\infty$, however
$$
K(r,0) = -\frac{\psi_{rr}}{\psi}(r,0) \rightarrow -1
$$
and
$$
K(r,\theta)  \rightarrow - \infty  \,\, ;\,\,\,\theta \neq 0 
$$
when $r \rightarrow +\infty$. Clearly 
$$
\left|\frac{\psi_{\theta}}{\psi}(r,\theta)\right|= r^{3}g'(\theta) = r^{3}(\theta/2)\frac{\sin(\theta/2)}{\theta/2}\cos(\theta/2)\leq Cr^{3}e^{-ar}\leq b
$$
in $C(a)$, therefore by Theorem (\ref{thm3}) the essential spectrum de $\mathbb{R}^{2}$ with such metric contains the interval $[1/4,+\infty)$.

\section{Acknowledgments}
This paper is part of the first
author's PhD thesis. L.M. acknowledges support from CAPES-Brazil and would like to thank the
Department of Mathematics at Universidade Federal do Cear\'a for the pleasant and productive
period during his PhD program at that institution. F.M. acknowledges support from CNPq-Brazil. 

\bibliographystyle{model3-num-names}

\begin{thebibliography}{00}



\bibitem{Davies} Davies, E. B., Spectral theory and differential operators, Cambridge University Press (1995).

\bibitem{Zhou} De Tang Zhou, Essential spectrum of the Laplacian on manifolds of nonnegative curvature, Int. Math. Res. Not. 5
(1994).

\bibitem{Donnelly} Donnelly, H., On the essential spectrum of a complete Riemannian manifold, Topology 20 (1981) 1-14.

\bibitem{Donn} Donnelly, H., Exhaustion functions and the spectrum of Riemannian manifolds, Indiana Univ. Math. J. 46 (2)
(1997) 505-527.

\bibitem{DoLi} Donnelly, H. and Li, P., Pure point Spectrum and Negative Curvature for Noncompact Manifolds. Duke Math. J. 46 (1979), 497-503.

\bibitem{Escobar} Escobar, Jos\'e F.; Freire, Alexandre, The spectrum of the Laplacian of manifolds of positive curvature, Duke Math. J.65 (1992), 1-21.

\bibitem{Kumura} Kumura, H., On the spectrum of the Laplacian on complete manifolds, J. Math. Soc. Japan 49 (1997), 1-14.

\bibitem{wang} Jiaping Wang, The spectrum of the Laplacian on a minifold of nonnegative Ricci curvature, Math. Res. Lett. 4 (4) (1997) 473-479.

\bibitem{Li} Li, J., Spectrum of the Laplacian on a complete Riemannian manifold with non-negative Ricci curvature which possesses a pole, J. Math. Soc. Japan 46 (1994), 213-216.

\bibitem{Chen} Zhi Hua Chen, Zhi Qin Lu, Essential spectrum of complete Riemannian manifolds, Sci. China Ser. A 35 (3) (1992), 276-282.
 
\bibitem{Zhiqin and Zhou} Zhiqin Lu, Detang Zhou, On the essential spectrum of complete non-compact
manifolds, Journal of functional analysis 260 (2011) 3283-3298.

\end{thebibliography}



\end{document}